\documentclass[12pt,verbatim]{amsart}
\usepackage{amsmath,amssymb,amsfonts}  
\usepackage{graphicx}
\usepackage{caption}
\usepackage{amsthm}
\usepackage{subfigure}
\usepackage{verbatim}
\usepackage{amscd}
\usepackage{ifthen}

\newtheorem{theorem}[equation]{Theorem}

\newtheorem{lemma}[equation]{Lemma}

\newtheorem{corollary}[equation]{Corollary}
\newtheorem{remark}[equation]{Remark}

\newtheorem{nonsec}[equation]{}

\newcommand{\Rn}{ {\mathbb{R}^n} }

\newcommand{\sh}{\,\textnormal{sh}}
\newcommand{\ch}{\,\textnormal{ch}}
\renewcommand{\th}{\,\textnormal{th}}

\numberwithin{equation}{section}

\title{A Gromov Hyperbolic metric and M\"obius transformations}
\author{Xiaoxue Xu} 
\address{School of Science, Zhejiang Sci-Tech University, Hangzhou 310018, China}  
\email{xiaoxue\_xu@126.com} 

\author{Gendi Wang*} 
\address{School of Science, Zhejiang Sci-Tech University, Hangzhou 310018, China}
\email{gendi.wang@zstu.edu.cn(*corresponding author)}

\author{Xiaohui Zhang} 
\address{School of Science, Zhejiang Sci-Tech University, Hangzhou 310018, China}
\email{xiaohui.zhang@zstu.edu.cn}

\begin{document}  

\newcounter{minutes}\setcounter{minutes}{\time}
\divide\time by 60
\newcounter{hours}\setcounter{hours}{\time}
\multiply\time by 60 \addtocounter{minutes}{-\time}
\def\thefootnote{}
\footnotetext{ {\tiny File:~\jobname.tex,
          printed: \number\year-\number\month-\number\day,
          \thehours.\ifnum\theminutes<10{0}\fi\theminutes }}
\makeatletter\def\thefootnote{\@arabic\c@footnote}\makeatother

\maketitle

\begin{abstract}
We compare a Gromov hyperbolic metric with the hyperbolic metric in the unit ball or in the upper half space,
and prove sharp comparison inequalities between the Gromov hyperbolic metric and some hyperbolic type metrics.
We also obtain several sharp distortion inequalities for the Gromov hyperbolic metric under some families of M\"{o}bius transformations.
\end{abstract}

{\small \sc Keywords.} {Gromov hyperbolic metric, hyperbolic metric, hyperbolic type metrics, M\"obius transformations }

{\small \sc 2010 Mathematics Subject Classification.} {30F45 (51M10)}

\section{Introduction}

It is well known that metrics play important roles in  geometric function theory.
One of the most important metrics is the hyperbolic metric in the unit ball or in the upper half space.
In addition to the classical hyperbolic metric, numerous hyperbolic type metrics are natural generalizations of the hyperbolic metric.
The most important property of the hyperbolic metric is its invariance under a group of M\"{o}bius transformations.
Examples of M\"obius invariant metrics also include the Seittenranta metric  \cite{s99},
the Apollonian metric \cite{b98}, and the M\"obius invariant Cassinian metric \cite{i19}.
In order to better understand these metrics, various estimates between the hyperbolic metric and hyperbolic type metrics are investigated \cite{avv,C,PH,P1,Ibragimov3,s99,Matti,z}.

The most used hyperbolic type metrics are the quasihyperbolic metric and the distance ratio metric \cite{gp, go}.
Whereas both metrics are not M\"obius invariant, then it is natural to study the quasi-invariance properties for these metrics.
Namely, it would be interesting to obtain the Lipschitz constants for these metrics under M\"obius transformations.
Indeed, Gehring, Palka, and Osgood have proved that the quasihyperbolic metric and the distance ratio metric are not changed by more than a constant $2$ under M\"obius transformations,
see \cite[Corollary 2.5]{gp} and \cite[proof of Theorem 4]{go}.
Several authors have also studied this topic for other hyperbolic type metrics in \cite{C,hvz,i19,klvw,MS2,MS,svw,wv,x}.

Recently,
Ibragimov introduced a new metric $ u_{Z} $ to hyperbolize the locally compact noncomplete metric space $ (Z,d) $ without changing its quasiconformal geometry which is defined as \cite{i0}
\begin{align*}
u_{Z}(x,y)
= 2\, \log \frac {d(x,y) + \max\{d(x,\partial Z),d(y,\partial Z)\}}
{\sqrt{d(x,\partial Z)\,d(y,\partial Z)}}\,,
\quad\quad x,y \in Z\,,
\end{align*}
where $ d(x,\partial Z) $ is the distance from the point $ x $ to the boundary of $Z$.
For a domain $ D \subsetneq \mathbb{R}^{n} $ equipped with the Euclidean metric,
we have \cite{MS2}
\begin{align*}
u_{D}(x,y)
= 2\, \log \frac {|x-y|+\max\{d(x),d(y)\}}
{\sqrt{d(x)\,d(y)}}\,,
\quad\quad x,y \in D,
\end{align*}
where $ d(x) $ denotes the Euclidean distance from $ x $ to the boundary of $D$.

Several authors have studied comparison inequalities between the Gromov hyperbolic metric and the hyperbolic metric
as well as some hyperbolic type metrics \cite{i0, MS2, z}.
Mohapatra and Sahoo also investigated quasi-invariance properties of the Gromov hyperbolic metric under quasiconformal mappings \cite{MS}.

In this paper, we continue the investigation on the Gromov hyperbolic metric to improve or complement some results in \cite{MS2}.
We further obtain sharp comparison inequalities between the Gromov hyperbolic metric and the hyperbolic metric, the distance ratio metric, and some other related hyperbolic type metrics such as the Seittenranta metric, the half-Apollonian metric and so on.
We also prove sharp distortion inequalities for the Gromov hyperbolic metric under some specific families of M\"obius transformations.

\section{Preliminaries}\label{section 2}
In this section, for readers' convenience, we collect the definitions and some basic properties of various hyperbolic type metrics.

\begin{nonsec}{\bf Hyperbolic metric.}
{\rm
The hyperbolic metrics $ \rho_{\mathbb{B}^{n}} $ and $ \rho_{\mathbb{H}^{n}} $
of the unit ball $ \mathbb{B}^{n} = \{ z \in \mathbb{R}^{n}: |z|<1 \} $ and
of the upper half space $ \mathbb{H}^{n} = \{x= (x_{1},\ldots,x_{n}) \in \mathbb{R}^{n}: x_{n}>0 \} $
are defined as follows.
By \cite[p.40]{AF}, for $ x,y \in \mathbb{B}^{n} $\,,
\begin{align*}
\sh \frac {\rho_{\mathbb{B}^{n}}(x,y)} {2}
=\frac {|x-y|} { \sqrt{(1-|x|^{2})(1-|y|^{2})}}\,,
\end{align*}
and hence \cite[(3)]{z}
\begin{align}\label{equ3}
\rho_{\mathbb{B}^{n}}(x,y)
=2 \log \frac {\sqrt{|x-y|^{2} + (1-|x|^{2})(1-|y|^{2})} + |x-y|} {\sqrt{(1-|x|^{2})(1-|y|^{2})}}\,.
\end{align}
By \cite[p.35]{AF}, for $ x,y \in \mathbb{H}^{n} $\,,
\begin{align*}
\ch \rho_{\mathbb{H}^{n}}(x,y)
=1 + \frac {|x-y|^{2}} {2 x_{n}y_{n}}\,,
\end{align*}
and hence \cite[(8)]{z}
\begin{align}\label{equ4}
\rho_{\mathbb{H}^{n}}(x,y)
=\log \left( 1 + \frac {|x-y|^{2} + \sqrt{|x-y|^{4} + 4x_{n}y_{n}|x-y|^{2}}} {2x_{n}y_{n}}\right)\,.
\end{align}
Two special formulas of the hyperbolic metric are frequently used \cite[(2.17),(2.6)]{Matti}\,:
\begin{align*}
\rho_{\mathbb{B}^{n}} (re_{1},se_{1})
=\log \left( \frac {1+s} {1-s} \cdot \frac {1-r} {1+r} \right)\,,
\quad {\rm for} \quad -1<r<s<1\,, \quad s>0\,,
\end{align*}
and
\begin{align*}
\rho_{\mathbb{H}^{n}} (re_{n},se_{n})
=\log\frac{s}{r}\,,
\quad {\rm for} \quad 0<r<s\,.
\end{align*}
}
\end{nonsec}

The following lemma shows the relation between the metric $u_{D}$ and the metric $\rho_{D}$ when $D\in\{\mathbb{B}^{n}\,, \mathbb{H}^{n}\}$\,.

\begin{theorem} \label{z1} \cite[Theorem 1, Theorem 2]{z}, \cite[Theorem 3.6]{MS2}
\begin{equation*}
\rho_{\mathbb{B}^{n}}(x,y)\le u_{\mathbb{B}^{n}}(x,y) \leq 3\, \rho_{\mathbb{B}^{n}}(x,y)\,,
\quad  for  \quad x, y \in\mathbb{B}^{n}\,.
\end{equation*}
\begin{equation*}
\rho_{\mathbb{H}^{n}}(x,y)\le u_{\mathbb{H}^{n}}(x,y) \leq 3\, \rho_{\mathbb{H}^{n}}(x,y)\,,
\quad  for  \quad x, y \in\mathbb{H}^{n}\,.
\end{equation*}
All the inequalities are sharp.
\end{theorem}

\begin{nonsec}{\bf Distance ratio metric.}
{\rm
In \cite[p.51]{go}, Gehring and Osgood introduced the distance ratio metric $ \tilde{j}_{D} $\,.
Let $D$ be a proper open subset  of $ \mathbb{R}^{n} $. For $ x,y \in D $,
\begin{align*}
\tilde{j}_{D}(x,y)
=\frac {1} {2} \, \log \left( 1+\frac {|x-y|} {d(x)} \right) \left( 1 + \frac {|x-y|} {d(y)} \right)\,.
\end{align*}

Vuorinen made some modification of the above definition and defined the metric $j_{D}$, still called the  distance ratio metric, as follows \cite[(2.34)]{Matti}:
\begin{align*}
j_{D}(x,y)
=\log \left( 1 + \frac {|x-y|} {\min\{d(x),d(y)\}} \right)\,.
\end{align*}

}
\end{nonsec}

The following lemma shows the relation between the distance ratio metric and the hyperbolic metric.

\begin{lemma}\label{c} \cite[Lemma 7.56]{avv}\label{le8}, \cite[Lemma 2.41(2)]{Matti}
\begin{equation*}
\frac12 \rho_{\mathbb{B}^{n}}(x,y)\le j_{\mathbb{B}^{n}}(x,y)
\leq\rho_{\mathbb{B}^{n}}(x,y)\,,
\quad  for  \quad x,y \in\mathbb{B}^{n}\,.
\end{equation*}
\begin{equation*}
\frac12 \rho_{\mathbb{H}^{n}}(x,y)\le j_{\mathbb{H}^{n}}(x,y)
\leq\rho_{\mathbb{H}^{n}}(x,y)\,,
\quad  for  \quad x,y \in\mathbb{H}^{n}\,.
\end{equation*}
\end{lemma}


\begin{theorem}\label{o3}\cite[Lemma 3.1]{MS2}
Let $D\subsetneq\mathbb{R}^{n}$ be arbitrary. Then
$$
2\,\tilde{j}_{D}(x,y) \leq u_{D}(x,y) \leq 4\,\tilde{j}_{D}(x,y)\,.
$$
The first inequality becomes equality when $d(x)=d(y)$.
\end{theorem}

\begin{theorem}\label{o4}\cite[Theorem 4.8]{MS2}
For $D\subsetneq\mathbb{R}^{n}$\,, we have
$$
j_{D}(x,y) \leq u_{D}(x,y) \leq 4\,j_{D}(x,y)\,.
$$
The first inequality is sharp.
\end{theorem}

\begin{nonsec}{\bf Seittenranta metric.}
{\rm
For an open subset $ D $ of $ \overline{\mathbb{R}}^{n} $ with card$ (\partial D) \geq 2 $ and for all
$ x,y \in D $\,,
the Seittenranta metric $ \delta_{D} $ is defined as \cite{s99}
\begin{align*}
\delta_{D}(x,y)
=\log \left( 1 + \sup_{p,q\in\partial D}|p,x,q,y| \right)\,,
\end{align*}
where
$$
|p,x,q,y|
=\frac {|p-q| \, |x-y|} {|p-x| \, |q-y|}
\quad {\rm with} \quad
\frac{|\infty-q|}{|\infty-x|}=1
$$
is the absolute ratio.

The most important property of the absolute ratio is its invariance under M\"obius transformations \cite[Theorem 3.2.7]{AF}.
It follows from the definitions that \cite[Remarks 3.2(3)]{s99}
$$
\delta_{\mathbb{R}^{n} \setminus \{\zeta\}} (x,y)= j_{\mathbb{R}^{n} \setminus \{\zeta\}} (x,y)
$$
for all $ \zeta \in \mathbb{R}^{n} $\,.
}
\end{nonsec}

The distance ratio metric and the Seittenranta metric are comparable as the following lemma shows.

\begin{lemma}\cite[Theorem 3.4]{s99}\label{js}
The inequalities
$$
j_{D}(x,y) \leq \delta_{D}(x,y)
\leq 2\, \tilde{j}_{D}(x,y)
\leq 2\, j_{D}(x,y)
$$
hold for every  open set $ D\subsetneq\mathbb{R}^{n} $\,.
\end{lemma}

\begin{nonsec}{\bf Apollonian metric.}
{\rm
For a proper open subset $ D $ of  $\overline{\mathbb{R}}^{n} $ and for all $ x,y \in D $\,,
the Apollonian metric $ \alpha_{D} $ is defined as \cite{b98}
\begin{align*}
\alpha_{D}(x,y)
=\sup_{p,q\in\partial D} \log|p,x,y,q|\,.
\end{align*}
}
\end{nonsec}

Note that $ \alpha_{D} $ is a pseudo-metric in $ D $\,.
It is, in fact, a metric if and only if $ \overline{\mathbb{R}}^{n} \setminus D $
is not contained in an $ (n-1) $-dimensional sphere in $ \overline{\mathbb{R}}^{n} $
\cite[Theorem 1.1]{b98}.
By \cite[Lemma 8.39]{Matti} and \cite[Example 3.2, Lemma 3.1]{b98}, we have
$$\delta_{D}(x,y)=\alpha_{D}(x,y)=\rho_{D}(x,y)$$
when $ D\in\{\mathbb{B}^{n}\,, \mathbb{H}^{n}\} $\,.
\medskip


The following lemma shows the relation between the metric $\alpha_{D}$ and the metric $j_{D}$.

\begin{lemma}\cite[Theorem 4.2]{s99}\label{alj}
Let $ D \subsetneq \mathbb{R}^{n} $ be a convex domain.
Then 
\begin{align*}
\alpha_{D}(x,y)
\leq j_{D}(x,y)\,.
\end{align*}
\end{lemma}

The following lemma shows the relation between the metric $\delta_{D}$ and the metric $\alpha_{D}$.

\begin{lemma}\cite[Theorem 3.11]{s99}\label{PAH}
Let $ D \subset \overline{\mathbb{R}}^{n} $ be an open set with {\rm card}$(\partial D) \geq2 $\,.
Then
\begin{align*}
\alpha_{D}(x,y)
\leq \delta_{D}(x,y)
\leq\log\left(e^{\alpha_{D}(x,y)}+2\right)
\leq\alpha_{D}(x,y)+\log3\,.
\end{align*}
The first two inequalities give the best possible bounds for $\delta_{D}$ expressed in terms of $\alpha_{D}$ only.
\end{lemma}

In \cite{PH}, H\"{a}st\"{o} and Lind\'{e}n gave another form of the Apollonian metric\,:
\begin{equation}\label{equ1}
\alpha_{D}(x,y)
= \sup_ { p \in \partial D } \log \frac {|p-y|} {|p-x|}
+ \sup_ { q \in \partial D } \log \frac {|q-x|} {|q-y|}\,.
\end{equation}
The half-Apollonian metric is defined by using one term in the right-hand side of (\ref{equ1}).

\begin{nonsec}{\bf Half-Apollonian metric.}
{\rm
For a proper open subset $ D $ of $ \mathbb{R}^{n} $ and for all $ x,y \in D $\,,
the half-Apollonian metric $ \eta_{D} $ is defined as \cite{PH}
\begin{align*}
\eta_{D}(x,y)
=\sup_{p\in\partial D} \left| \log \frac {|x-p|} {|y-p|} \right|\,.
\end{align*}
}
\end{nonsec}

Note that $ \eta_{D} $ is also a pseudo-metric in $ D $\,,
and a proper metric whenever
$ \mathbb{R}^{n}\setminus D $ is not a subset of a hyperplane \cite[Theorem 1.2]{PH}.
By \cite[Lemma 2.2 (i)]{b98}, we have
$$
\alpha_{D}(x,y) = \eta_{D}(x,y)
$$
when  $ D = \overline{\mathbb{R}}^{n} \setminus \{\zeta,\infty\}$
for any $ \zeta \in \mathbb{R}^{n} $\,.

\medskip

The following lemma shows the relation between the metric $\eta_{D}$ and the metric $\alpha_{D}$.

\begin{lemma} \cite[Theorem 2.1]{PH}\label{lemma3}
Let $ D \subsetneq \mathbb{R}^{n} $ be a domain.
Then the double inequality
\begin{equation*}
\frac{1}{2} \, \alpha_{D}(x,y) \leq \eta_{D}(x,y) \leq \alpha_{D}(x,y)
\end{equation*}
holds for all $ x,y \in D $\,.
Both inequalities are sharp.
\end{lemma}

\begin{nonsec}{\bf Cassinian metric.}
{\rm
For a proper subdomain $ D $ of $ \overline{\mathbb{R}}^{n} $ and for all $ x,y \in D $\,,
the Cassinian metric $ c_{D} $ is defined as \cite{i1}
\begin{align*}
c_{D}(x,y)
=\sup_{p\in\partial D} \frac {|x-y|} {|x-p||y-p|}\,.
\end{align*}
}
\end{nonsec}

\begin{nonsec}{\bf Triangular ratio metric.}
{\rm
For a proper subdomain $D$ of $\mathbb{R}^{n} $ and for all $ x,y \in D $\,,
the triangular ratio metric $ s_{D} $ is defined as \cite{C}
\begin{align*}
s_{D}(x,y)
=\sup_{p\in\partial D} \frac {|x-y|} {|x-p|+|y-p|}\,. 
\end{align*}
}
\end{nonsec}

The following lemma shows the relation between the metric $s_{D}$  and the metric $j_{D}$.

\begin{lemma}\label{equ14}\cite[Lemma 2.1]{hvz}
Let D be a proper subdomain of $\mathbb{R}^{n}$\,. Then
$$
\th\frac{j_{D}(x,y)}{2}
\leq s_{D}(x,y)
\leq\frac{e^{j_{D}(x,y)}-1}{2}\,.
$$
\end{lemma}

\bigskip

\section{The metric $u_D$ and the hyperbolic metric }

We devote this section to improving the right-hand side of inequalities (\ref{equ2}) and show the analogue result
in the upper half space.

\begin{theorem}\cite[Theorem 3.5]{MS2}
For all $ x,y \in \mathbb{B}^{n} $, we have
\begin{align}\label{equ2}
\rho_{\mathbb{B}^{n}}(x,y) - 2 \,\log2
\leq u_{\mathbb{B}^{n}}(x,y)
\leq 2\, \rho_{\mathbb{B}^{n}}(x,y) + 2\,\log2\,.
\end{align}
\end{theorem}

\medskip

\begin{theorem}\label{thm3.2}
For all $ x,y \in \mathbb{B}^{n} $\,, we have
\begin{align}\label{eqn:thm3.2}
u_{\mathbb{B}^n}(x,y)
\leq \rho_{\mathbb{B}^{n}}(x,y) + 2\,\log2\,,
\end{align}
and the inequality is sharp.
\end{theorem}
\begin{proof}
Without loss of generality,
we may assume that $ |x| \leq  |y|<1 $\,.
By (\ref{equ3}),
it suffices to prove that
\begin{align*}
\frac {|x-y|+1-|x|} {\sqrt{(1-|x|)(1-|y|)}}
\leq 2 \, \frac { |x-y| + \sqrt{|x-y|^{2} + (1-|x|^{2})(1-|y|^{2})}} {\sqrt{(1-|x|^{2})(1-|y|^{2})}}\,,
\end{align*}
which is equivalent to
\begin{align*}
\left(|x-y|+1-|x|\right)\sqrt{(1+|x|)(1+|y|)}
\leq 2 \, \left( |x-y| + \sqrt{|x-y|^{2}+(1-|x|^{2})(1-|y|^{2})} \right)\,.
\end{align*}
The above inequality follows from
\begin{align*}
\sqrt{(1+|x|)(1+|y|)}<2
\end{align*}
and
\begin{align*}
(1-|x|)^{2}
&\leq(1-|x||y|)^{2}
=(|x|-|y|)^{2}+(1-|x|^{2})(1-|y|^{2})\\
&\leq|x-y|^{2}+(1-|x|^{2})(1-|y|^{2})
\,.
\end{align*}
This proves the desired inequality.

For the sharpness, we set $x=se_1$ and $y=te_1$ with $0<s<t<1$. By (\ref{equ3}), we have
\begin{align}
\lim\limits_{t\to 1^{-}}\left(u_{\mathbb{B}^n}(x,y)-\rho_{\mathbb{B}^{n}}(x,y)\right)
={}&\lim\limits_{t\to 1^{-}}2\log\frac{(t-s+1-s)\sqrt{(1+s)(1+t)}}{t-s+\sqrt{(t-s)^2+(1-s^2)(1-t^2)}}\nonumber\\
={}&2\log\sqrt{2(1+s)}.\label{eqn:u-rho}
\end{align}
For arbitrary $\epsilon>0$, there exists $0<s_0<1$ such that $2\log\sqrt{2(1+s_0)}>2\log2-\epsilon/2$.
It follows from \eqref{eqn:u-rho} that there exists a number $t_0$ with $s_0<t_0<1$ such that
$$2\log\frac{(t_0-s_0+1-s_0)\sqrt{(1+s_0)(1+t_0)}}{t_0-s_0+\sqrt{(t_0-s_0)^2+(1-s_0^2)(1-t_0^2)}}>2\log\sqrt{2(1+s_0)}-\epsilon/2>2\log2-\epsilon.$$
Hence for arbitrary $\epsilon>0$, there exist $x=s_0e_1$ and $y=t_0e_1$ such that
$$u_{\mathbb{B}^n}(x,y)-\rho_{\mathbb{B}^{n}}(x,y)>2\log2-\epsilon,$$
which implies the sharpness of the inequality \eqref{eqn:thm3.2}.
\end{proof}

\medskip

\begin{theorem}\label{thm3.5}
For all $ x,y \in \mathbb{H}^{n} $\,,
we have
\begin{align*}
u_{\mathbb{H}^n}(x,y)
\leq \rho_{\mathbb{H}^{n}}(x,y) + 2\,\log2\,,
\end{align*}
and 
the inequality is sharp.
\end{theorem}
\begin{proof}
Without loss of generality,
we may assume that $ 0 < y_{n} \leq x_{n} $\,.
By (\ref{equ4}), it suffices to prove that
\begin{align*}
\left( \frac {|x-y|+x_{n}} {\sqrt{x_{n}y_{n}}} \right)^{2}
\leq 4\, \left( 1+\frac {|x-y|^{2}+|x-y|\sqrt{|x-y|^{2}+4\,x_{n}y_{n}}} {2\, x_{n}y_{n}} \right)\,,
\end{align*}
which is equivalent to
\begin{align*}
{x_{n}}^{2} + 2\, |x-y|x_{n}
\leq 4\, x_{n}y_{n} + |x-y|^{2} + 2\, |x-y| \sqrt{|x-y|^{2} + 4\, x_{n}y_{n}}\,.
\end{align*}
The above inequality follows from
\begin{align*}
{x_{n}}^{2}+2|x-y|x_{n}
&\leq (|x-y| + y_{n})^{2}+2|x-y|x_{n}\\
&\leq 4\, x_{n}y_{n} + |x-y|^{2} + 2\, |x-y|(x_{n}+y_{n})
\end{align*}
and
\begin{align*}
(x_{n} + y_{n})^{2}
\leq |x-y|^{2} + 4\, x_{n}y_{n}\,.
\end{align*}

For the sharpness, 
we choose $ x=e_{n} $ and $ y=te_{n} $ with $ 0<t<1\,.$
Then
\begin{align*}
\lim_{t\rightarrow 0^{+}} \left( u_{\mathbb{H}^{n}}(x,y) - \rho_{\mathbb{H}^{n}}(x,y) \right)
=\lim_{t\rightarrow 0^{+}} \left( 2\, \log\frac{2-t}{\sqrt{t}} - \log\frac{1}{t} \right)
=2\log2\,.
\end{align*}
Thus completes the proof.
\end{proof}

\bigskip

\section{The metric $u_D$ and the distance ratio metric }

Mohapatra and Sahoo \cite{MS2} compared the metric $u_{D}$ with the distance ratio metrics $\tilde{j}_{D}$ and $j_{D}$.
By Theorem \ref{o3} and Lemma \ref{js},
we have
\begin{align}\label{equ11}
j_{D}(x,y)
\leq 2\, \tilde{j}_{D}(x,y)
\leq u_{D}(x,y)
\leq 4\, \tilde{j}_{D}(x,y)
\leq 4\, j_{D}(x,y)\,,
\end{align}
and the first two inequalities give the best possible bounds for the metric $u_{D}$ in terms of the metrics $\tilde{j}_{D}$ and  ${j}_{D}$\,, see Theorem \ref{o3} and Theorem \ref{o4}.

\medskip


In this section, we will refine the last two inequalities in \eqref{equ11}. Specifically, the following Theorem \ref{jg} and Corollary \ref{jg4} show that the constant $4$ in inequalities \eqref{equ11} can be improved to $3$
and $3$ is the best possible.

\medskip

\begin{theorem}\label{jg}
For all $ x,y \in D \subsetneq \mathbb{R}^{n}$\,,
we have
\begin{align}\label{equ7}
u_{D}(x,y)\leq 3\, \tilde{j}_{D}(x,y)\,,
\end{align}
and the inequality is sharp.
\end{theorem}
\begin{proof}
Without loss of generality, we may assume that $ d(y) \leq  d(x) $\,.
To show our claim,
it suffices to prove that
\begin{align*}
\left( \frac {|x-y|+d(x)} {\sqrt{d(x)d(y)}} \right) ^{4}
\leq \left( 1 + \frac {|x-y|} {d(x)} \right) ^{3} \left( 1 + \frac {|x-y|} {d(y)} \right) ^{3}\,,
\end{align*}
which is equivalent to
\begin{align*}
|x-y| d(x) + d^{2}(x)
\leq \frac {|x-y|^{3}} {d(y)} + 3\, |x-y|^{2} + 3\, |x-y|d(y) + d^{2}(y)\,.
\end{align*}
The above inequality follows from
\begin{align*}
|x-y| d(x) \leq |x-y|(|x-y|+d(y)) = |x-y|^{2}+|x-y|d(y)
\end{align*}
and
\begin{align*}
d^{2}(x) \leq (|x-y|+d(y))^{2} = |x-y|^{2} + 2\, |x-y|d(y) + d^{2}(y)\,.
\end{align*}

For the sharpness,
we consider the domain $ D = \mathbb{R}^{n} \setminus \{e_{1}\} $\,.
Let $ x = -\,y = t\, e_{1} $ with $ 0 < t < 1 $\,.
Then
\begin{align*}
\lim_{t\rightarrow0^{+}} \frac{u_{D}(x,y)}{\tilde{j}_{D}(x,y)}
=\lim_{t\rightarrow0^{+}} \frac{2\log\frac{1+3t}{\sqrt{1-t^{2}}}}
{\frac{1}{2} \log \left( 1 +  \frac{2t}{1-t} \right) \left( 1 + \frac{2t}{1+t} \right)}
=3\,.
\end{align*}
Thus completes the proof.
\end{proof}

\medskip

\begin{corollary}\label{jg4}
For all $ x,y \in D \subsetneq \mathbb{R}^{n}$\,,
we have
\begin{align}\label{equ8}
u_{D}(x,y)
\leq 3\,j_{D}(x,y)\,,
\end{align}
and the inequality is sharp.
\end{corollary}
\begin{proof}
The inequality follows from Theorem \ref{jg} and Lemma \ref{js}\,.

For the sharpness,
we consider the domain $ D = \mathbb{R}^{n} \setminus \{e_{1}\} $\,.
Let $ x = -\, y = t\, e_{1} $ with $ 0 < t < 1 $\,.
Then
\begin{align*}
\lim_{t\rightarrow0^{+}} \frac{u_{D}(x,y)}{j_{D}(x,y)}
=\lim_{t\rightarrow0^{+}} \frac{2 \log\frac{1+3t}{\sqrt{1-t^{2}}}}{\log\left(1+\frac{2t}{1-t}\right)}
=3\,.
\end{align*}
Thus completes the proof.
\end{proof}

\begin{remark}
{\rm
By Lemma \ref{c}, it is obvious that Corollary \ref{jg4} improves
the upper bounds of the metric $u_D$ in \cite[Theorem 1, Theorem 2]{z} and \cite[Theorem 3.6]{MS2}, see Theorem \ref{z1}.
Moreover, the inequalities \eqref{equ7} and \eqref{equ8} hold in arbitrary proper subdomains of $\mathbb{R}^{n}$
}
\end{remark}

\begin{theorem}\label{hg}
For all $ x,y \in D \subsetneq \mathbb{R}^{n}$\,,
we have
\begin{align*}
u_{D}(x,y)
\leq 2\, \tilde{j}_{D}(x,y) + \log2\,,
\end{align*}
and 
the inequality is sharp.
\end{theorem}
\begin{proof}
Without loss of generality, we may assume that $ d(y) \leq d(x) \,.$
To show our claim,
it suffices to prove that
\begin{align*}
\frac {\left(|x-y|+d(x)\right)^{2}} {d(x)d(y)}
\leq 2\, \left( 1+\frac{|x-y|}{d(x)} \right) \left(1+\frac{|x-y|}{d(y)} \right)\,,
\end{align*}
or, equivalently
\begin{align*}
d(x) \leq |x-y| + 2\, d(y)\,,
\end{align*}
which is true by the triangle inequality.

For the sharpness, 
we consider the domain $ D = \mathbb{R}^{n} \setminus \{e_{1}\} $\,.
Let $ x = -\, y = t\,e_{1} $ with $ 0 < t < 1 $\,.
Then
\begin{align*}
 & \lim_{t\rightarrow 1^{-}} \left( u_{D}(x,y) - 2\, \tilde{j}_{D}(x,y) \right)\\
=& \lim_{t\rightarrow 1^{-}} \left( 2\log\frac{1+3t}{\sqrt{1-t^{2}}}
-\log \left( 1+\frac{2t}{1-t} \right) \left( 1+\frac{2t}{1+t} \right) \right)\\
=& \log2\,.
\end{align*}
Thus completes the proof.
\end{proof}

Theorem \ref{hg} and Lemma \ref{js} together yield the following corollary.

\begin{corollary}
For all $ x,y \in D \subsetneq \mathbb{R}^{n}$\,,
we have
\begin{align*}
u_{D}(x,y)
\leq 2\, j_{D}(x,y) + \log2\,.
\end{align*}
\end{corollary}

\section{The metric $u_D$ and other related metrics }

In this section,
we compare the metric $u_{D}$ with the metrics $ \delta_{D}$\,, $\eta_{D}$\,, $\alpha_{D}$\,,
$c_{D} $\,, and $ s_{D} $\,, respectively.

\begin{theorem}\label{equ13}\cite[Corollary 5.4]{MS2}
For all $x,y\in D\subsetneq\mathbb{R}^{n}$\,, we have
$$
\frac{1}{2}\,\delta_{D}(x,y)
\leq u_{D}(x,y)
\leq 4\, \delta_{D}(x,y)\,.
$$
\end{theorem}

The following theorem is an improvement of Theorem \ref{equ13} and of importance in studying the distortion property of the metric $u_D$ under M\"obius transformations.

\begin{theorem}\label{jg6}
For all $ x,y \in D \subsetneq \mathbb{R}^{n}$\,,
we have
\begin{align*}
\delta_{D}(x,y)
\leq u_{D}(x,y)
\leq 3\, \delta_{D}(x,y)\,,
\end{align*}
and both inequalities are sharp.
\end{theorem}
\begin{proof}
The inequalities follow from Theorem \ref{o3}, Lemma \ref{js} and Corollary \ref{jg4}.

For the sharpness of the left-hand side of the inequalities,
we consider the domain $D=\mathbb{B}^{n} $\,.
Let $ x = -\, y = t\, e_{1} $ with $ 0<t<1 $\,.
Then
$$
u_{D}(x,y) = \delta_{D}(x,y) = 2\, \log\frac{1+t}{1-t}\,.
$$

For the sharpness of the right-hand side of the inequalities,
we consider the domain $D = \mathbb{R}^{n} \setminus \{e_{1}\} $\,,
then $ \delta_{D}(x,y) = j_{D}(x,y)\,. $
By Corollary \ref{jg4}, the constant $3$ is the best possible.
\end{proof}

\begin{theorem}\label{jg8}\cite[Lemma 5.14]{MS2}
Let $ D\subsetneq\mathbb{R}^{n}$ and $x,y\in D$\,. Then
$$
\eta_{D}(x,y)
\leq u_{D}(x,y)
\leq 4 \log\left(2+e^{\eta_{D}(x,y)}\right)\,.
$$
\end{theorem}

\begin{theorem}\label{jg7}
For all $ x,y \in D \subsetneq \mathbb{R}^{n}$\,,
we have
\begin{equation*}
\eta_{D}(x,y) \leq u_{D}(x,y) \leq 2\, \eta_{D}(x,y) + 2\log3\,.
\end{equation*}
The constant $1$ in the left-hand side and the constant $ 2\log3 $ in the right-hand side of the inequalities are the best possible.
\end{theorem}
\begin{proof}
The left-hand side of the inequalities is the fact of Theorem \ref{jg8}.

For the sharpness of the constant $1$\,, we consider the domain
$ D = \mathbb{R}^{n} \setminus \{e_{1}\} $\,.
Let $ x=0 $ and $ y=t\,e_{1} $ with $0<t<1 $\,.
Then
\begin{align*}
\eta_{D}(x,y)
=\log \frac{1}{1-t}
\quad {\rm and} \quad
u_{D}(x,y)
= 2\, \log \frac {1+t} {\sqrt{1-t}} \,.
\end{align*}
Now we see that
\begin{align*}
\lim_{t\rightarrow1^{-}} \frac {u_{D}(x,y)} {\eta_{D}(x,y)}
&=\lim_{t\rightarrow1^{-}} \frac{2\, \log \frac {1+t} {\sqrt{1-t}}} {\log \frac{1}{1-t}}
=1\,.
\end{align*}

To prove the right-hand side of the inequalities,
we assume that $ d(x) \leq  d(y) $\,.
Choose $ z \in \partial D $ such that $ |x-z| = d(x) $\,.
This implies $ |x-z| \leq |y-z| $ and $ |x-y| \leq 2\,|y-z|\,.$
Then
\begin{align*}
u_{D}(x,y)
&=2\, \log \left( \frac {|x-y|+d(y)} {\sqrt{d(x)d(y)}} \right)
\leq 2\, \log \left( \frac {|x-y|+d(y)} {|x-z|} \right)\\
&\leq 2\, \log \left( \frac {3\,|y-z|} {|x-z|} \right)
\leq 2\, \sup_{w \in \partial D} \log \left( \frac {3\,|y-w|} {|x-w|} \right)\\
&= 2\, \eta_{D}(x,y)+2\, \log3\,.
\end{align*}

For the sharpness of the constant $2\log 3$\,, 
let $ D = \mathbb{R}^{n} \setminus\{0\} $ and $ y = -\, x\,.$ Then $u_{D}(x,y)=2\log 3$ and $\eta_{D}(x,y)=0$.

Thus completes the proof.
\end{proof}

\medskip

\begin{remark}
{\rm
Let
$$
f(t)=4 \, \log\left(2+e^{t}\right) - 2\, t-2\,\log3\,,
$$
where $t=\eta_{D}(x,y)\in[0,+\infty)$\,.
By differentiation, we have
$$
f'(t)=\frac{2\,e^{t}-4}{2+e^{t}}\,,
$$
which is negative on $(0,\log2)$ and positive on $(\log 2,+\infty)$.
Hence, we have
$$
f_{\min}(t)=f(\log2)=6\,\log2-2\,\log3>0\,.
$$
Therefore, the right-hand side of the inequalities in Theorem \ref{jg7} is better than that in Theorem \ref{jg8}.
}
\end{remark}

\begin{theorem}
For all $ x,y \in D \subsetneq \mathbb{R}^{n}$\,,
we have
\begin{align*}
\alpha_{D}(x,y)
\leq u_{D}(x,y)
\leq 2\,\alpha_{D}(x,y) + 2\, \log3 \,.
\end{align*}
The constant $1$ in the left-hand side and the constant $ 2\log3 $ in the right-hand side of the inequalities are the best possible.
\end{theorem}
\begin{proof}
The left-hand side of the inequalities follows from Theorem \ref{jg6} and Lemma \ref{PAH}.

For the sharpness of the constant $1$\,,
we consider the domain $ D = \mathbb{R}^{n} \setminus \{e_{1}\} $\,,
then $ \alpha_{D}(x,y) = \eta_{D}(x,y) \,.$
The result follows from Theorem \ref{jg7}.

The right-hand side of the inequalities follows from Theorem \ref{jg7} and Lemma \ref{lemma3}.

For the sharpness of the constant $2\log 3$\,,
we consider the domain $ D = \mathbb{R}^{n} \setminus\{0\} $\,, then $\alpha_{D}(x,y)=\eta_{D}(x,y)$\,.
The result follows from Theorem \ref{jg7}.
\end{proof}

\medskip

\begin{theorem}
Let $ D \subsetneq \mathbb{R}^{n} $ be a convex domain.
Then for all $x,y \in D$\,,
\begin{align*}
u_{D}(x,y)
\leq 3\,\alpha_{D}(x,y)\,,
\end{align*}
and the inequality is sharp.
\end{theorem}
\begin{proof}
The inequality follows from Corollary \ref{jg4} and Lemma \ref{alj}.

By Theorem \ref{z1}, the constant $ 3 $ is the best possible since $\alpha_D(x,y)=\rho_D(x,y)$ for $ D = \mathbb{B}^{n}$.
\end{proof}

\medskip

\begin{theorem}\label{equ12}
For all $ x,y \in D \subsetneq \mathbb{R}^{n}$\,,
we have
\begin{align*}
u_{D}(x,y)
\geq 2\, \log(1+r\,c_{D}(x,y))\,,
\end{align*}
where $ r=\min\{d(x),d(y)\} $\,.
The inequality is sharp.
\end{theorem}
\begin{proof}

The inequality follows from \cite[Theorem 4]{x}  and \cite[Theorem 4.5]{MS2}.

For the sharpness,
we consider the punctured space
$ D_{p}=\mathbb{R}^{n} \setminus\{p\} $\,.
Let $ x,y \in D_{p} $ with $ |x-p| = |y-p| $\,.
It is clear that
$$
u_{D}(x,y)
= 2\,\log \left( 1+\frac{|x-y|}{|x-p|} \right)
= 2\,\log \left( 1+|x-p| \, c_{D}(x,y) \right)\,.
$$
Hence the inequality is sharp.
\end{proof}

\begin{remark}
{\rm
By \cite[Corollary 5.6]{MS2}, we have
$$
u_{D}(x,y)\geq\frac{1}{2}\log\left(1+R\,c_{D}(x,y)\right)\,,
$$
where $R=\max\{d(x),d(y)\}$\,.
It is easy to see that the result in Theorem \ref{equ12} is better than that in \cite[Corollary 5.6]{MS2}
when $r\geq\frac{R}{4}$\,.
}
\end{remark}

\begin{theorem}
For all $ x,y \in D \subsetneq \mathbb{R}^{n}$\,,
we have
\begin{align*}
(2 \log3)\,s_{D}(x,y) \leq u_{D}(x,y)\leq 3 \, \log\frac{1+s_{D}(x,y)}{1-s_{D}(x,y)}\,,
\end{align*}
and both inequalities are sharp.
\end{theorem}
\begin{proof}
The left-hand side of the inequalities is the fact of \cite[Corollary 5.10]{MS2}.

For the sharpness of the left-hand side of the inequalities, we consider the domain $ D = \mathbb{R}^{n} \setminus\{0\} $\,.
Setting $y=-x$,
then $ s_{D}(x,y) = 1 $ and $ u_{D}(x,y) = 2\,\log3 $\,.


The right-hand side of the inequalities follows from Corollary \ref{jg4} and Lemma \ref{equ14}.

For the sharpness of the right-hand side of the inequalities,
we consider the domain $ D = \mathbb{R}^{n} \setminus \{e_{1}\} $\,.
Let $ y = -\, x = -t\, e_{1} $ with $ 0 < t < 1 $\,.
Then
\begin{align*}
u_{D}(x,y)=2\, \log\frac{1+3t}{\sqrt{1-t^{2}}}
\quad {\rm and} \quad
s_{D}(x,y)=t\,.
\end{align*}
Moreover,
\begin{align*}
\lim_{t\rightarrow0^{+}} \frac{u_{D}(x,y)} {\log\frac{1+s_{D}(x,y)}{1-s_{D}(x,y)}}
=\lim_{t\rightarrow0^{+}}  \frac{2\,\log\frac{1+3t}{\sqrt{1-t^{2}}} }{\log\frac{1+t}{1-t}}
=3\,.
\end{align*}
Thus completes the proof.
\end{proof}

\section{the metric $u_D$ and M\"{o}bius transformations }

In this section, we study quasi-invariance properties for the metric $u_D$ under M\"{o}bius transformations.
We first give the distortion inequalities for the metric $u_{D}$ under M\"{o}bius transformations in arbitrary domains $D\subsetneq\mathbb{R}^{n}$.

\begin{theorem}\label{M}
Let $ D $ and $ D' $ be proper subdomains of $ \mathbb{R}^{n} $ and
$ f: \overline{{\mathbb R}}^n \to \overline{{\mathbb R}}^n $ be a M\"{o}bius transformation with
$ f D = D' $. Then for $ x,y \in D $, we have
\begin{align}\label{45}
\frac{1}{3}\, u_{D}(x,y)
\leq u_{D'}(f(x),f(y))
\leq 3\,u_{D}(x,y)\,.
\end{align}
\end{theorem}
\begin{proof}
The proof follows from the inequalities in Theorem \ref{jg6} and M\"obius invariance of the metric $\delta_D$.
\end{proof}

Now we discuss the sharpness of inequalities \eqref{45} in some specific domains.

\begin{theorem}
Let $ f $ be a M\"{o}bius transformation with $ f\mathbb{B}^{n} = \mathbb{B}^{n} $\,.
Then for all $ x,y \in \mathbb{B}^{n} $\,,
we have
\begin{align*}
\frac{1}{3}\, u_{\mathbb{B}^{n}}(x,y)
\leq u_{\mathbb{B}^{n}}(f(x),f(y))
\leq 3\, u_{\mathbb{B}^{n}}(x,y)\,,
\end{align*}
and both inequalities are sharp.
\end{theorem}
\begin{proof}
The double inequality is clear by Theorem \ref{M}.

For the sharpness of the right-hand side of the inequalities, we consider
$$
f(z)=a^{\ast}+\frac{r^{2}(z-a^{\ast})}{|z-a^{\ast}|^{2}}\,,
$$
where $a=te_{1}$  $(\frac{1}{2} < t <1 )$\,,
$a^{\ast}=\frac{a}{|a|^{2}}$\,,
$r=\sqrt{|a^{\ast}|^{2}-1}$\,.

Putting $ x = -\, y = (1-t)\, e_{1} $\,.
Then
$$
f(x)
=\frac{2 \, t-1}{1-t+t^{2}}e_{1}
\quad {\rm and} \quad
f(y)
=\frac{1}{1+t-t^{2}}e_{1}\,.
$$
Moreover,
\begin{align*}
\lim_{t\rightarrow1^{-}} \frac{u_{\mathbb{B}^{n}} \left(f(x),f(y) \right)}{u_{\mathbb{B}^{n}}(x,y)}
=\lim_{t\rightarrow1^{-}} \frac{\log\frac{4+t-5 \, t^{2}+t^{3}}{\sqrt{t(2-t)(1-t+t^{2})(1+t-t^{2})}}}
{\log\frac{2-t}{t}}=3\,.
\end{align*}

Thus the constant $3$ is attained.
The sharpness of the left-hand side of the inequalities can be seen by considering the inverse of $f$ and hence the constant $\frac{1}{3}$ is also the best possible.

\end{proof}


\begin{theorem}
Let $ f $ be a M\"{o}bius transformation with $ f\mathbb{H}^{n} = \mathbb{B}^{n} $\,.
Then for all $ x,y \in \mathbb{H}^{n} $\,,
we have
\begin{align*}
\frac{1}{3}\, u_{\mathbb{H}^{n}}(x,y)
\leq u_{\mathbb{B}^{n}}(f(x),f(y))
\leq 3\, u_{\mathbb{H}^{n}}(x,y)\,,
\end{align*}
and the left-hand side of the inequalities is sharp.
\end{theorem}
\begin{proof}
The double inequality is clear by Theorem \ref{M}.

For the sharpness of the left-hand side of the inequalities, we consider
\begin{equation*}
f(z)
=-\, e_n + \frac{2(z+e_n)}{|z+e_n|^2}\,.
\end{equation*}

Putting $ x = t\, e_{n} $ and $ y = \frac{1}{t}\, e_{n}$ with $ t>1 $\,.
Then
\begin{align*}
f(x)
=-\, \frac{t-1}{t+1}\, e_{n}
\quad {\rm and} \quad
f(y)
=\frac{t-1}{t+1}\, e_{n}\,.
\end{align*}
Moreover,
\begin{align*}
\lim_{t\rightarrow1^{+}} \frac{u_{\mathbb{B}^{n}}(f(x),f(y))}{u_{\mathbb{H}^{n}}(x,y)}
=\lim_{t\rightarrow1^{+}} \frac{2\, \log t}{2\, \log\left(2t-\frac{1}{t}\right)}
=\frac{1}{3}\,.
\end{align*}
Thus completes the proof.
\end{proof}


Next we give another type distortion inequalities for the metric $u_{D}$ under M\"{o}bius transformations in some specific domains.

\begin{theorem}\label{thm6}
Let $ f$ be a M\"{o}bius transformation with $ f\mathbb{H}^{n} = \mathbb{H}^{n} $\,.
Then for all $ x,y \in \mathbb{H}^{n} $\,,
we have
\begin{align*}
u_{\mathbb{H}^{n}}(x,y) - 2\, \log2
\leq u_{\mathbb{H}^{n}}(f(x),f(y))
\leq u_{\mathbb{H}^{n}}(x,y) + 2\, \log2\,,
\end{align*}
and both inequalities are sharp.
\end{theorem}
\begin{proof}
By Theorem \ref{thm3.5} and Theorem \ref{z1}, we obtain
\begin{align*}
u_{\mathbb{H}^{n}}(f(x),f(y))
&\leq \rho_{\mathbb{H}^{n}}(f(x),f(y)) + 2\, \log2
= \rho_{\mathbb{H}^{n}}(x,y) + 2\, \log2\\
&\leq u_{\mathbb{H}^{n}}(x,y) + 2\, \log2
\end{align*}
and
\begin{align*}
u_{\mathbb{H}^{n}}(f(x),f(y))
&\geq \rho_{\mathbb{H}^{n}}(f(x),f(y))
=\rho_{\mathbb{H}^{n}}(x,y)\\
&\geq u_{\mathbb{H}^{n}}(x,y) - 2\, \log2\,.
\end{align*}

For the sharpness of the left-hand side of the inequalities,
we consider
$$
f(z)=\frac{z}{|z|^{2}}\,.
$$

Putting $ x = e_{1} + t \, e_{n} $ and
$ y = e_{1} + \frac{1}{t} \, e_{n} $ with $ 0<t<1 $\,.
Then
$$
f(x)= \frac {1} {1+t^{2}} \, e_{1} + \frac {t} {1+t^{2}} \, e_{n}
\quad\quad {\rm and} \quad\quad
f(y)= \frac {t^{2}} {1 + t^{2}} \, e_{1} + \frac {t} {1+t^{2}} \, e_{n}\,.
$$
Moreover,
$$
\lim_{t\rightarrow 0^{+}} \left( u_{\mathbb{H}^{n}}(x,y) - u_{\mathbb{H}^{n}}(f(x),f(y))  \right)
=\lim_{t\rightarrow 0^{+}} \left( 2\log\left(\frac{2}{t}-t\right)
-2\log\left(1+\frac{1-t^{2}}{t}\right) \right)
=2\log2\,.
$$

Thus the left-hand side of the inequalities is sharp.
The sharpness of the right-hand side of the inequalities can be seen by considering the inverse of $f$\,.

\end{proof}

\begin{theorem}
Let $ f $ be a M\"{o}bius transformation with $ f\mathbb{H}^{n} = \mathbb{B}^{n} $\,.
Then for all $ x,y \in \mathbb{H}^{n} $\,,
we have
\begin{align*}
u_{\mathbb{H}^{n}}(x,y) - 2\, \log2
\leq u_{\mathbb{B}^{n}}(f(x),f(y))
\leq u_{\mathbb{H}^{n}}(x,y) + 2\,\log2\,,
\end{align*}
and the left-hand side of the inequalities is sharp.
\end{theorem}
\begin{proof}
By Theorem \ref{thm3.2}, Theorem \ref{thm3.5} and Theorem \ref{z1}, we obtain
\begin{align*}
u_{\mathbb{B}^{n}}(f(x),f(y))
&\leq \rho_{\mathbb{B}^{n}}(f(x),f(y)) + 2\, \log2
= \rho_{\mathbb{H}^{n}}(x,y) + 2\, \log2\\
&\leq u_{\mathbb{H}^{n}}(x,y) + 2\, \log2
\end{align*}
and
\begin{align*}
u_{\mathbb{B}^{n}}(f(x),f(y))
&\geq \rho_{\mathbb{B}^{n}}(f(x),f(y))
=\rho_{\mathbb{H}^{n}}(x,y)\\
&\geq u_{\mathbb{H}^{n}}(x,y)-2\, \log2\,.
\end{align*}

For the sharpness of the left-hand side of the inequalities, we consider
\begin{equation*}
f(z)
=-\, e_n + \frac{2(z+e_n)}{|z+e_n|^2}\,.
\end{equation*}

Putting $ x = t\, e_{n} $ and $ y = \frac{1}{t}\, e_{n}$ with $ t>1 $\,.
Then
\begin{align*}
f(x)
=-\, \frac{t-1}{t+1}\, e_{n}
\quad {\rm and} \quad
f(y)
=\frac{t-1}{t+1}\, e_{n}\,.
\end{align*}
Moreover,
\begin{align*}
\lim_{t\rightarrow \infty} \left( u_{\mathbb{H}^{n}}(x,y) - u_{\mathbb{B}^{n}}(f(x),f(y)) \right)
=\lim_{t\rightarrow \infty} \left( 2\, \log\left(2t-\frac{1}{t}\right) - 2\, \log t  \right)
=2\, \log2 \,.
\end{align*}
Thus completes the proof.
\end{proof}

\begin{theorem}
Let $ f $ be a M\"{o}bius transformation with $ f\mathbb{B}^{n} = \mathbb{B}^{n} $\,.
Then for all $ x,y \in \mathbb{B}^{n} $\,,
we have
\begin{align*}
u_{\mathbb{B}^{n}}(x,y) - 2\, \log2
\leq u_{\mathbb{B}^{n}}(f(x),f(y))
\leq u_{\mathbb{B}^{n}}(x,y) + 2\,\log2\,.
\end{align*}
\end{theorem}
\begin{proof}
By Theorem \ref{thm3.2} and Theorem \ref{z1}, we obtain
\begin{align*}
u_{\mathbb{B}^{n}}(f(x),f(y))
&\leq \rho_{\mathbb{B}^{n}}(f(x),f(y)) + 2\, \log2
= \rho_{\mathbb{B}^{n}}(x,y) + 2\, \log2\\
&\leq u_{\mathbb{B}^{n}}(x,y) + 2\, \log2
\end{align*}
and
\begin{align*}
u_{\mathbb{B}^{n}}(f(x),f(y))
&\geq \rho_{\mathbb{B}^{n}}(f(x),f(y))
=\rho_{\mathbb{B}^{n}}(x,y)\\
&\geq u_{\mathbb{B}^{n}}(x,y)-2\, \log2\,.
\end{align*}
Thus completes the proof.
\end{proof}

\subsection*{Acknowledgments}
This research was partly supported by National Natural Science Foundation of China (NNSFC) under Grant No.11771400 and No.11601485,
and Science Foundation of Zhejiang Sci-Tech University (ZSTU) under Grant No.16062023\,-Y.


\end{document}